\documentclass[12pt]{amsart}
\usepackage{amsmath,amsfonts,amssymb}

\newtheorem{theorem}{Theorem}[section]
\newtheorem{lemma}{Lemma}[section]
\newtheorem{corollary}{Corollary}[section]

\usepackage[english]{babel}
\def \R {\mathbb {R}}
\def \C {\mathbb {C}}

\def \p {\partial}
\def \O {\Omega}

\newcommand{\tr}{\hbox{\rm trace\,}}

\begin{document}
\title{ \sc Integral formulas for a class of curvature PDE's \\  and applications to isoperimetric inequalities \\ and to symmetry problems }
\author{{ Vittorio Martino}}
\thanks{Dipartimento di Matematica, Universit\`a di Bologna,
piazza di Porta S.Donato 5, 40127 Bologna, Italy. E-mail address:
{\tt{martino@dm.unibo.it}}}
\author{{Annamaria Montanari}}
\thanks{Dipartimento di Matematica, Universit\`a di Bologna,
piazza di Porta S.Donato 5, 40127 Bologna, Italy. E-mail address:
{\tt{montanar@dm.unibo.it}}}
\date{}

\begin{abstract}{We prove  integral formulas for closed hypersurfaces in
$\C^{n+1},$ which furnish a relation between elementary symmetric
functions in the eigenvalues of the complex Hessian matrix of the
defining function and the Levi curvatures of the hypersurface. Then
we follow the Reilly approach to prove an isoperimetric inequality.
As an application, we obtain the ``Soap Bubble Theorem'' for
star-shaped domains with positive and constant Levi curvatures
bounding the classical mean curvature from above.}
\end{abstract}
\maketitle
\markboth{Vittorio Martino \& Annamaria
Montanari}{Integral formulas for a class of curvature PDE's}
\section{Introduction} The study of surfaces in the
Euclidean space with either constant Gauss curvature or constant
mean curvature received in the past a great amount of attention. In
$1899$ Liebmann \cite{LI} proved that the spheres are the only
compact surfaces in $\R^3$ with constant Gauss curvature. In $1952$
S\"uss \cite{S} extended the Liebmann result showing that a compact
convex hypersurface in the Euclidean space must be a sphere,
provided that for some $j$ the $j-$th elementary symmetric
polynomial in the principal curvatures is constant. In $1954$ Hsiung
\cite{HS} proved that the ``convexity'' assumption can be relaxed to
the ``star-shapedness'' one. The proofs of the above results are
based on the Minkowski formulae. A breakthrough for this sort of
problems was made by Alexandrov \cite{AL1} in $1956,$ who proved
that the sphere is the only compact hypersurface embedded into
 the Euclidean space with constant mean curvature. Alexandrov method is completely
 different from the Liebmann-S\"uss method, and is based on the moving plane technique,
 on the interior maximum principle for elliptic equations and on the boundary maximum
 principle of Hopf type for uniformly elliptic equations.
 In 1978 Reilly \cite{RE2} obtained another proof of the Alexandrov theorem combining the Minkowski formulae
 with some new elegant arguments.

It is well known that the Levi curvatures  contain less geometric
information than Euclidean curvatures, because the Levi form is only
a part of the second fundamental form. However, in a joint paper
with Lanconelli \cite{ML} we wrote the $j$-th Levi curvature for
real hypersurfaces in $\C^{n+1}$ in terms of elementary symmetric
functions of the eigenvalues of the normalized Levi form,
 and we proved a strong comparison principle, leading to symmetry theorems
for domains with constant curvatures.

Precisely:
\begin{theorem}\label{MontLan} Let $D\subseteq \C^{n+1}$ be a
strictly $j$-pseudoconvex domain with connected boundary,
 $1\leq j\leq n.$ Let $B_R(z_0)\subseteq D$ be a tangent sphere to
  $\p D$ at some point $p\in \p D.$ If $K_{\p D}^{(j)}(z)$ is the  $j$-th Levi curvature of $\p D$
at $z\in \p D$ and
\[
K_{\p D}^{(j)}(z)\geq {1}/{R}^j, \quad \forall z \in \p D,
\]
then
 $D=B_R(z_0).$
\end{theorem}
In \cite{ML} we prove that if   $\O$ is a bounded domain of
 $\C^{n+1},$ with boundary a real hypersurface of class $C^2,$
 then the  $j$-th Levi curvature of $\p \O$
at $z=(z_1,\dots,z_{n+1})\in \p \O$ writes in term of defining
function $f$ of $\O=\{f(z)<0\}$ as
\begin{equation}\label{k1}
  K_{\p \O}^{(j)}(z)=-\frac{1}{\begin{pmatrix}
 n\\
 j
\end{pmatrix}}\frac{1}{|\p f|^{j+2}}\sum_{1\leq i_1<\dots<i_{j+1}\leq n+1}
\Delta_{(i_1,\cdots,i_{j+1})}(f)
\end{equation}
for all $j=1,\dots,n$, where $|\p
f|=\sqrt{\sum_{j=1}^{n+1}|f_j|^2}$
\begin{equation}\label{k2}
\Delta_{(i_1,\cdots,i_{j+1})}(f)=\det\left(
\begin{array}{llll}
 0 & f_{\bar i_1} & \ldots &  f_{\overline{i}_{j+1}}\\
 f_{i_1}  & f_{i_1, \bar i_1} & \ldots & f_{i_1, \overline{i}_{j+1}} \\
  \vdots & \vdots & \ddots & \vdots \\
 f_{i_{j+1}}  & f_{i_{j+1}, \bar i_1} & \ldots & f_{i_{j+1}, \overline{i}_{j+1}}
\end{array}%
\right)
\end{equation}
and $f_j=\frac{\p f}{\p z_j},$ $f_{\bar j}=\overline{f_j},$ $f_{j
\bar \ell}=\frac{\p^2 f}{\p z_j \p \bar z_\ell}.$

Theorem \ref{MontLan} suggested the following question: are
spheres the unique compact hypersurfaces with constant Levi
curvatures? Klingenberg in \cite{K} gave a first positive answer
to this problem  by showing that a compact and strictly
pseudoconvex real hypersurface
 $M\subset \C^{n+1}$ is isometric to a sphere, provided that $M$ has constant horizontal
 mean curvature and the CR structure $T_{1,0}(M)$ is parallel in $T^{1,0}(\C^{n+1})$.
Later on in \cite{MaM} we relaxed Klingerberg conditions and we
proved that if  the characteristic direction is a geodesic, then
Alexandrov Theorem holds for hypersurfaces with positive constant
Levi mean curvature.

The problem of characterizing hypersurfaces with constant Levi
curvature has been studied by many authors. Hounie and Lanconelli
in \cite{HL} showed the result for Reinhardt domain of $\C^2,$
i.e. for domains $D$ such that if $ (z_1,z_2)\in D$ then $ (e^{i
\theta_1} z_1, e^{i \theta_2} z_2 )\in D$ for all real $\theta_1,
\theta_2 .$ Under this hypothesis, in a neighborhood of a point,
there is a defining function $F$ only depending on the radius
$r_1=|z_1|,$ $r_2=|z_2|,$ $F(r_1,r_2)=f(r_2^2)-r_1^2$ with $f$ the
solution of the ODE
\begin{equation}\label{ODE}
  sff''=sf'^2-k(f+sf'^2)^{3/2}-ff'
\end{equation}
 {Alexandrov Theorem follows from uniqueness of the solution
of \eqref{ODE}.}  Their technique has then been used in \cite{HL2}
to prove an Alexandrov Theorem for Reinhardt domains in $\C^{n+1}$
with spherical symmetry for every $n.$ Then in \cite{MM} Monti and
Morbidelli proved a Darboux -type theorem for $n\geq 2$: the unique
Levi umbelical hypersurfaces in $\C^{n+1}$ with all constant Levi
curvatures are spheres or cylinders.

In this paper we prove some integral formulas for compact
hypersurfaces, which furnish a relation between elementary symmetric
functions in the eigenvalues of the complex Hessian matrix of the
defining function and the Levi curvatures of the hypersurface. Then
we follow the Reilly approach to prove the following

\begin{theorem}[\sc Isoperimetric estimates]\label{stima}
Let $\O$ be a bounded domain of $\C^{n+1}$ with boundary a real
hypersurface of class $C^\infty.$ If $K^{(j)}_{\p \O}$ is non
negative at every point of $\p \O$ then
\begin{equation}
  \label{sigmacurvature3}
 \int_{\p \O}\left( \frac{1}{K_{\p
  \O}^{(j)}(x)}\right)^{1/j}d\sigma(x)\geq{2(n+1)}|\O|
\end{equation}
where $|\O|$ is the  Lebesgue measure of $\O.$ If $K^{(j)}_{\p
\O}$ is constant, then the equality holds  in
\eqref{sigmacurvature3} if and only if $\O$ is a ball of radius
$\left( \frac{1}{K_{\p
  \O}^{(j)}}\right)^{1/j}.$\par

\end{theorem}

The isoperimetric estimates \eqref{sigmacurvature3} show that the
Levi curvatures, which have been first introduced in analogy with
Euclidean curvatures by Slodkowski and Tomassini in \cite{Tom} and
\cite{STom}, have a deep geometric meaning, because they contain
information about the measure of a bounded domain.

Let us remark that there are non spherical sets which satisfy the
equality in \eqref{sigmacurvature3} (see \eqref{asterisco}). Thus,
the class of sets which satisfy the equality in
\eqref{sigmacurvature3}  is larger than the class  of sets which
satisfy the equality in the classical isoperimetric inequality and
in the Alexandrov Fenchel inequalities for quermassintegrals ( see
\cite{F}, \cite{G} and \cite{T}).

On the other side, if $K^{(j)}_{\p \O}$ is constant, then
\begin{equation}
  \label{sigmacurvaturecost}
 \left(K_{\p
  \O}^{(j)}\right)^{1/j}\leq\frac{|\p \O|}{{2(n+1)}|\O|}
\end{equation}
and an Alexandrov type theorem holds for star-shaped domains whose
classical mean curvature is bounded from above by a constant
$j$-Levi curvature. Precisly

\begin{corollary}[\sc An Alexandrov type Theorem]\label{alex}
Let  $\O\subset \C^{n+1}$ be a bounded star-shaped domain with
boundary a smooth real hypersurface. If the $j$-Levi curvature is a
positive constant $K^{(j)}$ at every point of $\p\O,$ then the
maximum of the mean curvature of $\p \O$ is bounded from below by
$(K^{(j)})^{(1/j)}.$ Moreover, if the mean curvature of $\p \O$ is
bounded from above by $(K^{(j)})^{(1/j)}$, then $\p \O$ is a sphere
and $\O$ is a ball.
\end{corollary}

Let us remark that, even if the Levi form is a part of the second
fundamental form, in general it is not possible to bound  from above
the Levi curvatures
 with the Euclidean ones. Indeed, it is very easy to build a cylinder
 in $\C^2$
 whose Levi curvature is $1/2$ while the mean curvature is $1/3.$

Our paper is organized as follows.
 In Section 2 we use the null Lagrangian property for elementary symmetric functions
  in the eigenvalues of the complex Hessian matrix and the classical divergence theorem to prove
  an integral formula
  \eqref{sigmacurvature} for a closed hypersurface in term of the
 $j$-th Levi curvature. In Section 3 we prove the isoperimetric
 estimates
 \eqref{sigmacurvature3}
 and we use the Minkowski formula
for the classical mean curvature  to conclude the proof of
Corollary \ref{alex}.
\par
\noindent {\bf Acknowledgement} We would like to thank Ermanno
Lanconelli and  Daniele Morbidelli for some useful discussions about
this argument.

\section{Null Lagrangian property for elementary symmetric functions
in the eigenvalues of the complex  Hessian matrix}

Given a Hermitian  matrix  $A,$ let   $\sigma_j(A)$ be the $j$-th
elementary symmetric function in the eigenvalues of $A$. Let us
recall that if $A$ is the  $(n+1)\times (n+1)$ matrix with
eigenvalues $\lambda_1,\dots, \lambda_{n+1}$ then
\[
\sigma_j(A)= \sum_{1\leq i_1<\dots<i_j\leq
n+1}\lambda_{i_1}\cdots\lambda_{i_j}.
\]
If we choose $A=[a_{\ell \bar k}]=\p\bar \p f$ the complex Hessian
matrix of a smooth function $f$ and we denote by $\frac{\p
\sigma_{j}(A)}{\p a_{\ell \bar k}}$ the partial derivative of the
function $\sigma_j$ with respect to the term of place ${\ell \bar
k}$ of the matrix $A$, then
\begin{equation}\label{null}
\sum_{\ell}\partial_\ell\left(\frac{\p \sigma_{j+1}(\p \bar \p
f)}{\p a_{\ell \bar k}}\right)=0, \forall k=1,\dots, n+1.
\end{equation}
(see \cite{R2} for a similar argument for the real Hessian
matrix).

Moreover, by using the notation \eqref{k2}, the following Lemma
holds
\begin{lemma}
  For every $f\in C^2$ and for every $j=1,\dots,n+1$
  \begin{equation}\label{omog2}
    \sum_{\ell, k=1}^{n+1}\frac{\p\sigma_{j+1}}{\p a_{\ell \bar k
    }}(\p\bar\p f)f_\ell f_{\bar k}=-\sum_{1\leq i_1<\dots<i_{j+1}\leq n+1}\Delta_{(i_1,
    \dots,
    i_{j+1})}(f)
  \end{equation}
\end{lemma}
\begin{proof}
  By explicitly writing $\sigma_{j+1}(\p\bar\p f)$ in the left hand side of \eqref{omog2} as $$\sigma_{j+1}(\p\bar\p f)=\sum_{1\leq i_1<\dots<i_{j+1}\leq n+1}\left|%
\begin{array}{ccc}f_{i_1, \bar i_1} & \ldots & f_{i_1, \overline{i}_{j+1}} \\
  \vdots & \ddots & \vdots \\
 f_{i_{j+1}, \bar i_1} & \ldots & f_{i_{j+1}, \overline{i}_{j+1}}
\end{array}%
\right|$$
with $|A|=\det A $ for every Hermitian matrix $A,$
we get
  \begin{equation}\label{omog3}
    \begin{split}
 \sum_{\ell, k=1}^{n+1}&\frac{\p\sigma_{j+1}}{\p a_{\ell \bar k
    }}(\p\bar\p f)f_\ell f_{\bar k}=
     \sum_{\ell, k=1}^{n+1}f_\ell f_{\bar k}\sum_{1\leq i_1<\dots<i_{j+1}\leq n+1}\frac{\p }{\p a_{\ell \bar k
    }}\left|%
\begin{array}{ccc}
 f_{i_1, \bar i_1} & \ldots & f_{i_1, \overline{i}_{j+1}} \\
  \vdots & \ddots & \vdots \\
 f_{i_{j+1}, \bar i_1} & \ldots & f_{i_{j+1}, \overline{i}_{j+1}}
\end{array}%
\right|\\
&=\sum_{1\leq i_1<\dots<i_{j+1}\leq n+1}\left(
 \sum_{\ell, k\in\{i_1,\dots,i_{j+1}\}}
 \frac{\p }{\p a_{\ell \bar k
    }}\left|%
\begin{array}{ccc}
 f_{i_1, \bar i_1} & \ldots & f_{i_1, \overline{i}_{j+1}} \\
  \vdots & \ddots & \vdots \\
 f_{i_{j+1}, \bar i_1} & \ldots & f_{i_{j+1}, \overline{i}_{j+1}}
\end{array}%
\right|f_\ell f_{\bar k}\right).\\
    \end{split}
  \end{equation}

On the other side, if we put $F(\p f,\bar \p f,\p \bar \p f)=
  -\Delta_{(i_1, \dots,
    i_{j+1})}(f)$ and we twice differentiate it with respect to
    $f_{i_\ell}$ and $f_{\overline{i}_{k}}$
    for every $\ell, k=1,\dots, j+1$, we get

    \[
\frac{\p F}{\p f_{i_\ell}}=(-1)^{\ell+1} \left|
\begin{array}{ccc}
  f_{\bar i_1} & \ldots &  f_{\overline{i}_{j+1}}\\
  f_{i_1, \bar i_1} & \ldots & f_{i_1, \overline{i}_{j+1}} \\
  \vdots & \ddots & \vdots \\
  f_{i_{\ell-1}, \bar i_1} & \ldots & f_{i_{\ell-1}, \overline{i}_{j+1}}\\
  f_{i_{\ell+1}, \bar i_1} & \ldots & f_{i_{\ell+1}, \overline{i}_{j+1}}\\
 \vdots & \ddots & \vdots\\
f_{i_{j+1}, \bar i_1} & \ldots & f_{i_{j+1}, \overline{i}_{j+1}}
\end{array}
 \right|
    \]
    \begin{equation}\label{dersec}
\begin{split}
\frac{\p^2 F}{\p f_{i_\ell}\p f_{\bar i_k}}=&(-1)^{\ell+k} \left|
\begin{array}{cccccc}
 f_{i_1, \bar i_1} & \ldots & f_{i_1, \overline{i}_{k-1}}& f_{i_1, \overline{i}_{k+1}}& \ldots & f_{i_1, \overline{i}_{j+1}}\\
   \vdots & \ddots & \vdots &  \vdots & \ddots & \vdots  \\
  f_{i_{\ell-1}, \bar i_1} & \ldots & f_{i_{\ell-1}, \overline{i}_{k-1}}&  f_{i_{\ell-1}, \overline{i}_{k+1}} & \ldots & f_{i_{\ell-1}, \overline{i}_{j+1}}\\
   f_{i_{\ell+1}, \bar i_1} & \ldots & f_{i_{\ell+1}, \overline{i}_{k-1}}&  f_{i_{\ell+1}, \overline{i}_{k+1}} & \ldots & f_{i_{\ell+1}, \overline{i}_{j+1}} \\
   \vdots & \ddots & \vdots &  \vdots & \ddots & \vdots  \\
 f_{i_{j+1}, \bar i_1} & \ldots & f_{i_{j+1}, \overline{i}_{k-1}} & f_{i_{j+1}, \overline{i}_{k+1}} & \ldots & f_{i_{j+1},
 \overline{i}_{j+1}}
\end{array}
 \right|
 \\
=&\frac{\p }{\p a_{\ell \bar k
    }}\left|%
\begin{array}{ccc}
 f_{i_1, \bar i_1} & \ldots & f_{i_1, \overline{i}_{j+1}} \\
  \vdots & \ddots & \vdots \\
 f_{i_{j+1}, \bar i_1} & \ldots & f_{i_{j+1}, \overline{i}_{j+1}}
\end{array}%
\right|\\
 \end{split}
    \end{equation}
    Moreover,
    \[
\begin{split}
  F(\p f, \bar \p f, \p \bar \p f)=& F(\p f, \bar \p f, \p \bar \p
  f)-F(0, \bar \p f, \p \bar \p f)\\
  =&\int_0^1\frac{d}{ds}F(s\p f, \bar \p f, \p \bar \p f)ds\\
   =&\int_0^1 \left(\sum_{\ell\in\{i_1,\dots,i_{j+1}\}}\frac{\p F}{\p f_\ell}(s\p f, \bar \p f, \p \bar \p
   f)f_\ell\right)ds\\
   =&\int_0^1ds \left(\sum_{\ell\in\{i_1,\dots,i_{j+1}\}}\frac{\p F}{\p f_\ell}(\p f, \bar \p f, \p \bar \p
   f)f_\ell\right)\\
=&\sum_{\ell\in\{i_1,\dots,i_{j+1}\}}\frac{\p F}{\p f_\ell}(\p f,
\bar \p f, \p \bar \p
   f)f_\ell
\end{split}
    \]
and by the same argument and by \eqref{dersec}
\begin{equation}
\label{omog4}
\begin{split}
 F(\p f, \bar \p f, \p \bar \p f)
=&\sum_{\ell,k\in\{i_1,\dots,i_{j+1}\}}\frac{\p^2 F}{\p f_\ell \p
\bar f_k}(\p f, \bar \p f, \p \bar \p
   f)f_\ell f_{\bar k}\\
=&\sum_{\ell,k\in\{i_1,\dots,i_{j+1}\}}\frac{\p }{\p a_{\ell \bar k
    }}\left|%
\begin{array}{ccc}
 f_{i_1, \bar i_1} & \ldots & f_{i_1, \overline{i}_{j+1}} \\
  \vdots & \ddots & \vdots \\
 f_{i_{j+1}, \bar i_1} & \ldots & f_{i_{j+1}, \overline{i}_{j+1}}
\end{array}%
\right|f_\ell f_{\bar k}.\\
\end{split}
\end{equation}
By substituting \eqref{omog4} in \eqref{omog3} we get \eqref{omog2}.
\end{proof}

By using the null lagrangian property for elementary symmetric
functions in the eigenvalues of the complex Hessian matrix  and the
classical divergence Theorem we get the following integral formulas
for closed hypersurfaces.

\begin{theorem}\label{Anna}
Let $\O$ be a bounded domain of $\C^{n+1}$ with boundary a real
hypersurface of class $C^2.$ For every defining function $f$ of
class $C^2$ of $\O=\{f(z)<0\}$ and for every $j=1,\dots,n$ we have
\begin{equation}
  \label{sigmacurvature}
 \int_\O \sigma_{j+1}(\p \bar \p f)dx={\begin{pmatrix}
 n+1\\
 j+1
\end{pmatrix}}\frac{1}{ 2(n+1)}\int_{\p \O}K_{\p
  \O}^{(j)}(z)|\p f|^{j+1}d\sigma(x),
\end{equation}
where  $K_{\p
  \O}^{(j)}$ is the $j$-th Levi curvature of $\p \O.$
\end{theorem}

\begin{proof}

Since $\sigma_j$ is a homogenous function of degree $j,$ i.e.
$\sigma_j(tA)=t^j\sigma_j(A)$ for every real $t$, we get
\begin{equation}\label{omog}
  \sigma_{j+1}(\p \bar \p f)=\frac{1}{j+1}\sum_{\ell, k=1}^{n+1}
\frac{\p \sigma_{j+1}(\p \bar \p f)}{\p a_{\ell \bar k}}f_{\ell \bar
k}
\end{equation}
Let us set $\nu_\ell =\frac{\p_\ell f}{|\p f|}$ and and let us
identify $z\in \C^{n+1}$ with $x\in \R^{2(n+1)},$ then by
\eqref{omog}, by \eqref{null}, by the classical divergence theorem
and by \eqref{omog2} we get
  \[
  \begin{split}
    \int_\O \sigma_{j+1}(\p \bar \p f)dx &=\frac{1}{j+1}
 \int_\O \sum_{\ell,k=1}^{n+1}\p_\ell\left(
\frac{\p \sigma_{j+1}}{\p a_{\ell \bar k}}(\p \bar \p f)f_{ \bar
k}\right)dx\\
&=\frac{1}{2(j+1)}\int_{\p \O} \sum_{\ell,k=1}^{n+1}\left(
\frac{\p \sigma_{j+1}}{\p a_{\ell \bar k}}(\p \bar \p f)f_{ \bar
k}\nu_\ell\right)d\sigma(x)
\\
&=\frac{1}{2(j+1)}\int_{\p \O} \sum_{\ell,k=1}^{n+1}\frac{\left(
\frac{\p \sigma_{j+1}}{\p a_{\ell \bar k}}(\p \bar \p f)f_{ \bar
k}f_\ell\right)}{|\p f|}d\sigma(x)
\\
&=-\frac{1}{2(j+1)}\int_{\p \O}\frac{ \sum_{1\leq
i_1<\dots<i_{j+1}\leq n+1} \Delta_{(i_1,\cdots,i_{j+1})}(f)}{|\p
f|}d\sigma(x)
\\
&= {\begin{pmatrix}
 n\\
 j
\end{pmatrix}}\frac{1}{ 2(j+1)}\int_{\p \O}K_{\p
  \O}^{(j)}(z)|\p f|^{j+1}d\sigma(x)\\
&= {\begin{pmatrix}
 n+1\\
 j+1
\end{pmatrix}}\frac{1}{ 2(n+1)}\int_{\p \O}K_{\p
  \O}^{(j)}(z)|\p f|^{j+1}d\sigma(x).
\end{split}
\]
\end{proof}
\section{Isoperimetric estimates and an Alexandrov type theorem}
In this section we use the integral formula \eqref{sigmacurvature}
to get an estimate of the $j$-th Levi curvature of a closed
hypersurface and to show Theorem \ref{stima}.

\begin{proof}[Proof of Theorem \ref{stima}]
If $ \int_{\p \O}\big(1/{K^{(j)}_{\p
\O}}\big)^{1/j}d\sigma(x)=+\infty $ then the inequality
\eqref{sigmacurvature3} is trivial. Thus, in the sequel we shall
assume $ \int_{\p \O}\big(1/{K^{(j)}_{\p
\O}}\big)^{1/j}d\sigma(x)<+\infty.$
 Let $f:\bar \O \rightarrow \R$ be the $C^2(\bar \O)$ solution of the
 Dirichlet problem
  \begin{equation}\label{Laplacian}
    \left\{%
\begin{array}{ll}
    \tr \p \bar \p f=1, & \hbox{in }\O; \\
    f=0, & \hbox{on}\,\p \O. \\
\end{array}%
\right.
  \end{equation}
Let us remark that $\tr \p \bar \p =\frac{1}{4}\Delta,$ with
$\Delta$ the Laplacian operator in $\R^{2n+2}.$ Hence, if $\p \O$ is
of class $C^{2,\alpha},$ the Dirichlet problem \eqref{Laplacian} has
a unique solution $f\in C^2(\bar \O).$ Let us recall that for every
$(n+1) \times (n+1)$ Hermitian matrix $A$ the Newton inequality
holds
\begin{equation}\label{Newton}
\sigma_j(A)\leq {\begin{pmatrix}
 n+1\\
 j
\end{pmatrix}}\left(\frac{\tr A}{n+1}\right)^j
\end{equation}
for all $j=2,\dots,n+1.$ Moreover,  in \eqref{Newton} equality holds
iff the matrix    $A$ is proportional to the identity matrix.

By applying \eqref{Newton} to the complex Hessian matrix of $f,$
with $f$ a solution of \eqref{Laplacian}, we get an estimate from
above of the left hand side of \eqref{sigmacurvature}
\begin{equation}\label{left}
\begin{split}
  \int_\O \sigma_{j+1}(\p \bar \p f)dx &\leq {\begin{pmatrix}
 n+1\\
 j+1
\end{pmatrix}}\frac{1}{(n+1)^{j+1}}\int_\O \left(\tr(\p \bar \p f)\right)^{j+1}dx\\
&= {\begin{pmatrix}
 n+1\\
 j+1
\end{pmatrix}}\frac{|\O|}{(n+1)^{j+1}}.
\end{split}\end{equation}
By applying again the diverge Theorem we get
\[
\int_{\p \O}|\p f|d\sigma(x)= \frac{1}{2}\int_{\p \O}\langle
\nabla f, N \rangle d\sigma(x)=\frac{1}{2}\int_{ \O}\Delta f
dx=2|\O|
\]
and by using the Cauchy-Schwarz inequality in the right hand side of
\eqref{sigmacurvature} we get
\begin{equation}\label{right}
\begin{split}
{\begin{pmatrix}
 n+1\\
 j+1
\end{pmatrix}}\frac{1}{2(n+1)}
 \int_{\p \O}K^{(j)}_{\p \O}|\p f|^{j+1}d\sigma(x)&\geq
\frac{{\begin{pmatrix}
 n+1\\
 j+1
\end{pmatrix}}\left(\int_{\p \O}|\p
 f|d\sigma(x)\right)^{j+1}}{{2(n+1)}\left(\int_{\p
 \O}\left(\frac{1}{K^{(j)}_{\p \O}}\right)^{1/j}d\sigma(x)\right)^{j}}\\
&= \frac{{\begin{pmatrix}
 n+1\\
 j+1
\end{pmatrix}}(2|\O|)^{j+1}}{{2(n+1)}\left(\int_{\p
 \O}\left(\frac{1}{K^{(j)}_{\p \O}}\right)^{1/j}d\sigma(x)\right)^{j}}.\\
\end{split}
\end{equation}
In \eqref{right} the equality   holds iff $|\p f|$ is proportional
to $\left(\frac{1}{K^{(j)}_{\p \O}}\right)^{1/j}.$ By the equality
\eqref{sigmacurvature} and by the inequalities \eqref{left} and
\eqref{right} we deduce
\[
\frac{(2|\O|)^{j}}{\left(\int_{\p
 \O}\left(\frac{1}{K^{(j)}_{\p \O}}\right)^{1/j}d\sigma(x)\right)^{j}}\leq
\frac{1}{(n+1)^{j}}
\]
and we get
\[
\int_{\p
 \O}\left(\frac{1}{K^{(j)}_{\p \O}}\right)^{1/j}d\sigma(x)\geq
2(n+1)|\O|.
\]
Moreover, in \eqref{left} the equality holds iff the complex Hessian
matrix of $f$ is proportional to the identity matrix.
 Since the defining function for $\O$ has been chosen such that   $\tr \p\bar \p
 f=1$,
 then it should be
  $\p\bar \p f=\frac{1}{n+1}I$ in $\bar
\O$ and by  \eqref{k1} and \eqref{k2}  we have
\begin{equation}\label{gradc}
  (K^{(j)}_{\p \O})^{1/j}=\frac{1}{(n+1)|\p f|}
\end{equation}
on $\p \O.$ The last equality does not seem enough to conclude that
$\O $ is a ball. Indeed, by the maximum principle for the Laplacian
operator, the function $f$ has an interior minimum at $z_0\in \O,$
and it is not restrictive to assume $z_0=0$ because the Levi
curvature equations are invariant with respect to translations.
Then, for every pluriharmonic function $h$ such that $h(0)=0$ and
$\p h(0)=0,$
\begin{equation}\label{asterisco}
  f(z)=f(0)+\frac{1}{n+1}|z|^2+h(z)
\end{equation}
is a defining function for $\O.$ For example, if in $\C^2$ we choose
\[
h(z_1,z_2)=\frac{{\rm Re}\, (z_1^2+z_2^2)}{4}
\]
then the set of zeros of the function $f$ in \eqref{asterisco} is
not a sphere.
 However, if $K^{(j)}_{\p \O}$ is constant
for some $j$, then by \eqref{gradc} $|\p f|$ should be constant on
$\p \O.$ It follows that the Dirichlet problem \eqref{Laplacian} is
over determinate  and by Serrin Theorem \cite{Se} we can conclude
that $\O$ is a ball and $\p \O$ is a sphere.
\end{proof}
Let  $H$ be the Euclidean mean curvature of  $\p \O.$ We recall that
the celebrated
 Minkowski formula (see for instance \cite{R}) asserts that
\begin{equation}\label{Vintegr}
\int_{\p \O} \,d\sigma=\int_{\p \O}H(x)\,  \langle N , x \rangle
d\sigma (x),
\end{equation} where $N$ is the outward unit normal.
\par
We are now ready to conclude the proof of Corollary \ref{alex}
\begin{proof}[Proof of Corollary \ref{alex}]
If $\O$ is star-shaped with respect to a point
 then by
\eqref{Vintegr} and by the divergence theorem we have
\begin{equation}\label{const}
\begin{split}
 |\p \O|=\int_{\p \O} \,d\sigma\leq \max_{\p \O}H\int_{\p \O}\langle N ,
x \rangle
 d\sigma (x)&=\max_{\p \O}H\int_\O\left(\sum_j\p_{x_j} x_j \right)dx\\
 &=2(n+1)\max_{\p \O}H|\O|.
\end{split} \end{equation}
By \eqref{const} we have
\[
\max_{\p \O}H\geq \frac{|\p \O|}{2(n+1)|\O|}
\]
and since $K^{(j)}$ is a positive constant, then by
\eqref{sigmacurvaturecost} we get
\[
(K^{(j)})^{1/j}\leq \frac{|\p \O|}{2(n+1)|\O|}\leq \max_{\p \O}H.
\]
Moreover, if $\max_{\p \O}H\leq (K^{(j)})^{1/j}$ then equality holds
in \eqref{sigmacurvaturecost} and by Theorem
 \ref{stima} we conclude that $\O$ is a ball.
\end{proof}

\bibliographystyle{alpha}

\end{document}